\documentclass[reqno]{amsart}

\usepackage{amsmath}
\usepackage{amssymb}
\usepackage{bbm}
\usepackage{mathrsfs}
\usepackage{float}

\newcommand{\br}{ }
\newcommand{\brr}{, }

\newtheorem{theo}{Theorem}[section]
\newtheorem{lem}[theo]{Lemma}
\newtheorem{fac}[theo]{Fact}
\newtheorem{prop}[theo]{Proposition}
\newtheorem{cons}[theo]{Conjectures}

\theoremstyle{definition}
\newtheorem{ttt}[theo]{}
\newtheorem{defi}[theo]{Definition}
\newtheorem{algo}[theo]{Algorithm}

\theoremstyle{remark}
\newtheorem{rem}[theo]{Remark}
\newtheorem{rems}[theo]{Remarks}
\newtheorem{nota}[theo]{Notation}
\newtheorem{ex}[theo]{Example}

\newtheorem{eo}[theo]{Experimental observation}


\def\pmod#1{\nobreak\ifinner\mkern8mu\else\mkern18mu\fi (\text{\rmfamily\upshape mod}\,\,#1)}

\newcommand{\liminv}{\mathop{\underleftarrow{\rm lim}}\limits}

\newcommand{\Pic}{\mathop{\text{\rm Pic}}\nolimits}
\newcommand{\Frob}{\mathop{\text{\rm Frob}}\nolimits}
\newcommand{\CH}{\mathop{\text{\rm CH}}\nolimits}
\newcommand{\tr}{\mathop{\text{\rm tr}}\nolimits}
\newcommand{\M}{\mathop{\text{\rm M}}\nolimits}
\newcommand{\N}{\mathop{\text{\rm N}}\nolimits}

\newcommand{\rk}{\mathop{\text{\rm rk}}\nolimits}
\newcommand{\cl}{\mathop{\text{\rm cl}}\nolimits}

\newcommand{\bbC}{{\mathbbm C}}
\newcommand{\bbF}{{\mathbbm F}}
\newcommand{\bbG}{{\mathbbm G}}
\newcommand{\bbN}{{\mathbbm N}}
\newcommand{\bbQ}{{\mathbbm Q}}
\newcommand{\bbZ}{{\mathbbm Z}}

\newcommand{\calO}{{\mathscr{O}}}

\newcommand{\frakp}{{\mathfrak{p}}}

\newcommand{\et}{{{\text{\rm {\'e}t}}}}
\newcommand{\tri}{{{\text{\rm tr}}}}
\newcommand{\id}{{{\text{\rm id}}}}

\newcommand{\Pb}{{\text{\bf P}}}
\newcommand{\Ab}{{\text{\bf A}}}

\newcounter{abc}
\newenvironment{abc}{\begin{list}{\rm \alph{abc}) }%
{\usecounter{abc} \leftmargin=0.0pt \labelsep=0.0pt %
\listparindent=0.0pt \labelwidth=0.0pt \parsep=\smallskipamount %
\itemsep=0.0pt \topsep=0.0pt \partopsep=\smallskipamount}}{\end{list}}

\newcounter{iii}
\newenvironment{iii}{\begin{list}{\rm \roman{iii}) }%
{\usecounter{iii} \leftmargin=0.0pt \labelsep=0.0pt %
\listparindent=0.0pt \labelwidth=0.0pt \parsep=\smallskipamount%
 \itemsep=0.0pt \topsep=0.0pt \partopsep=\smallskipamount}}{\end{list}}

\def\eop{\ifmmode\rule[-2pt]{0pt}{1pt}\ifinner\tag*{$\square$}\else\eqno{\square}\fi\else\vadjust{}\hfill$\square$\fi}

\title[Point counting on
$K3$~surfaces
and applications]{Point counting on
\boldmath{$K3$}~surfaces
and an application \\concerning real and complex multiplication}

\author{Andreas-Stephan Elsenhans}

\address{Institut f\"ur Mathematik \\ Warburger Stra\ss e 100 \\ D-33098 Paderborn, Germany}
\email{Stephan.Elsenhans@math.upb.de}
\urladdr{http://www.staff.uni-bayreuth.de/~bt270951/}

\author{J\"org Jahnel}

\address{\mbox{Department Mathematik\\ Univ.\ \!Siegen\\ \!Walter-Flex-Str.\ \!3\\ D-57068 \!Siegen\\ \!Germany}}
\email{jahnel@mathematik.uni-siegen.de}
\urladdr{http://www.uni-math.gwdg.de/jahnel}

\date{Mai~17,~2016.}
  
\keywords{point counting, Harvey's $p$-adic method, Weil polynomial, $K3$ surface, real and complex multiplication}

\subjclass[2010]{14J28 (primary), 11G15, 11Y16, 11M38 (secondary)}

\begin{document}

\begin{abstract}
We report on our project to find explicit examples of
$K3$~surfaces
having real or complex multiplication. Our strategy is to search through the arithmetic consequences of RM and CM. In~order to do this, an efficient method is needed for point counting on surfaces defined over finite~fields. For~this, we describe algorithms that are
\mbox{$p$-adic}
in~nature.
\end{abstract}

\maketitle

\section{Introduction}

\looseness-1
Let~$X$
be a quasi-projective variety over a finite
field~$\bbF_{\!q}$
of
characteristic~$p > 0$.
One~of the most natural questions concerning the arithmetic
of~$X$
is certainly to pinpoint the number of points that
$X$~has
over its base
field~$\bbF_{\!q}$.
For~example, let
$X := V(f) \subset \Pb^N$
be given
for~$f \in \bbF_{\!q}[T_0,\ldots,T_N]$
a homogeneous polynomial. Then~this means just to count the number of solutions in
$\bbF_{\!q}^{N+1}$
of the equation
$f(T_0,\ldots,T_N) = 0$,
ignoring about
$(0,\ldots,0)$
and up to~scaling.

More generally, one asks for the sequence
$(\#X(\bbF_{\!q^i}))_{i\ge1}$
of integers associated
with~$X$
that is composed of the numbers of points
on~$X$
that are defined over the extension
fields~$\bbF_{\!q^i}$.
This~sequence is of practical as well as theoretical interest. For~example, it is known since the days when E.\ Artin wrote his Ph.D.~thesis~\cite[\S22]{Ar} that it is wise to form the generating~function
\begin{equation}
\label{zeta}
Z_X(t) := \exp\left(\sum_{i=1}^\infty \#X(\bbF_{\!q^i}) \frac{t^i}i\right) ,
\end{equation}
which is called the zeta function of the
variety~$X$.
In~fact,
$Z_X$
is always a rational~function. According~to \cite[Exp.~XV, \S3, n$^\circ$2 and~3]{SGA5}, it is the alternating product of the Weil
poly\-no\-mi\-als~$\chi_i$
associated
with~$X$,
the characteristic polynomials of the operation of the geometric Frobenius
$\Frob$
on the
\mbox{$l$-adic}
cohomology vector spaces
$H_{c,\et}^i(X_{\overline\bbF_{\!q}}, \bbQ_l)$
with compact~support,
\begin{equation}
\label{ratzeta}
Z_X(t) = \!\!\prod_{j=0}^{2\dim X}\!\!\!\! \det \big(1-t\Frob | H_{c,\et}^j(X_{\overline\bbF_{\!q}}, \bbQ_l)\big)^{(-1)^{j+1}} =: \frac{\chi_1(t)\chi_3(t)\cdots\chi_{2\dim X-1}(t)}{\chi_0(t)\chi_2(t)\chi_4(t)\cdots\chi_{2\dim X}(t)} \, .
\end{equation}
Moreover,~at least when
$X$
is proper and smooth, the Weil polynomials
$\chi_j$
have remarkable properties. For~example, as was conjectured by A.~Weil and proven by P.~Deligne in \cite[Th\'eo\-r\`eme~(1.6)]{DeW1} and \cite[Corollaire~(3.3.9)]{DeW2}, every complex root
of~$\chi_j$
is of absolute value
$q^{j/2}$.
In~particular, no cancellations occur in formula~(\ref{ratzeta}).

Assume~now that the geometry of the variety
$X_{\overline\bbF_{\!q}}$
is well understood. By~which we mean that at least all its
\mbox{$l$-adic}
Betti numbers are~known. This~is usually the case in~practice, as the~varieties considered are curves of known genus, abelian varieties of known dimension, non-singular complete intersections of known multidegree,
$K3$~surfaces,~etc.
From~the algorithmic point of view, formula~(\ref{ratzeta}) then shows that it is easy to calculate
$\#X(\bbF_{\!q^i})$
for some very
large~$i$,
provided that
$\#X(\bbF_{\!q}), \#X(\bbF_{\!q^2}), \ldots, \#X(\bbF_{\!q^k})$
have been found. Thereby,~the
bound~$k$
only depends on the geometry
of~$X_{\overline\bbF_{\!q}}$.
Indeed,~once
$k$~is
properly chosen, (\ref{zeta}) and~(\ref{ratzeta}) allow to compute all coefficients of the
polynomials~$\chi_j$
from the numbers of points known. Then~to pinpoint
$\#X(\bbF_{\!q^i})$
is just the calculation of the
\mbox{$i$-th}
power series coefficient of
$\sum_j (-1)^{j+1} \log\chi_j$.

Unfortunately,~for most nontrivial types of varieties, the
bound~$k$
turns out to be too large, such that it is not feasible to count
$\#X(\bbF_{\!q}), \#X(\bbF_{\!q^2}), \ldots, \#X(\bbF_{\!q^k})$
in a naive way, even for medium-sized values
of~$q$.
For~instance, when
$X$~is
a
$K3$~surface,
things are relatively~simple. In~fact,
$\chi_1(t) = \chi_3(t) = 1$,
$\chi_0(t) = 1-t$,
$\chi_4(t) = 1-q^2t$,
and only
$\chi_2$
varies,
which is of 
degree~$22$.
Nevertheless,~one may only hope that
$k=10$
is~sufficient. It~is known that there are cases where
$k$~has
to be chosen larger~\cite[Table~6]{EJants9}.\smallskip

{\em Our motivation. Real and complex multiplication for
$K3$~surfaces.}
Our motivation comes from projective varieties
$X$
that are defined over a number
field~$K$.
As~was first noticed by R.~van Luijk, from point counts on the reductions
$X_{\frakp_1}$,
$X_{\frakp_2}$
modulo two primes~\mbox{\cite{vL05,vL07,EJants8,EJ12}} of good reduction or sometimes only one~\cite{EJ11}, one may determine the geometric N\'eron-Severi rank
of~$X$.
This~applies well even to varieties of general type~\cite{EK}, but a nontrivial case is provided already by
$K3$~surfaces.

The very basic idea behind van Luijk's method is that, in the case of a
$K3$~surface,
$\smash{\rk\Pic X_{\overline\bbF_{\!\frakp}}}$
is always even and one has
$\smash{\rk\Pic X_{\overline\bbF_{\!\frakp}} \geq \rk\Pic X_{\overline{K}}}$.
Moreover,~one hopes to find a prime such that
$\smash{\rk\Pic X_{\overline\bbF_{\!\frakp}}}$
is actually equal to
$\rk\Pic X_{\overline{K}}$
or~$\rk\Pic X_{\overline{K}} + 1$.
As~was observed by F.~Charles~\cite{Ch}, the existence of such primes is related to whether
$X$
has real multiplication or~not. More~precisely, existence is provided unless
$X$
has real multiplication and
$(22 - \rk\Pic X_{\overline{K}})/[E:\bbQ]$
is odd, for
$E$
the endomorphism~field.

We~say that a
$K3$~surface~$X$
has real or complex multiplication when the endomorphism algebra of the transcendental part
$T \subset H^2(X(\bbC), \bbQ)$,
considered as a pure
\mbox{$\bbQ$-Hodge}
structure, is strictly larger
than~$\bbQ$.
It~is then either a totally real field or a CM~field \cite[Theorem~1.6.a) and Theorem~1.5.1]{Za83}. Starting~from this definition, it seems, however, hard to conclude anything for concrete~examples. One~may deduce~\cite{vGe} at least that having RM or CM is a property of positive codimension in the analytic moduli space of
$K3$~surfaces
of fixed Picard rank (unless the Picard rank is~20, in which case every
$K3$~surface
has~CM). Thus,~having RM or CM should be thought of as being exceptional.

In~particular, it seems that one does not find examples of RM or CM~surfaces
just by~accident. Our~approach is therefore to perform our searches through the {\em arithmetic consequences\/} of RM and CM. These~include that, for all
primes~$\frakp$
that are (at least partially) inert
in~$E$,
the reduction
$X_\frakp$
is non-ordinary. I.e.,~that one has
$\#X_\frakp(\bbF_{\!\frakp}) \equiv 1 \pmod p$
for
$p$
the prime number
below~$\frakp$.
Moreover,~the transcendental factors of the corresponding Weil polynomials are all partially split over the same field, the endomorphism
field~$E$.
Cf.~the top of Section~\ref{RMCM} for a more precise~statement.
In~order to detect the very few surfaces showing such a behaviour within a large family, efficient methods for point counting are asked~for.\smallskip

{\em Harvey's
\mbox{$p$-adic}
method for counting points.}
The most efficient methods for counting points on a variety over a finite field that are known today are
\mbox{$p$-adic}
in nature. One~of them was originally developed by K.~Kedlaya in~\cite{Ke}. Kedlaya's point of view was cohomological. This~means, he computed the characteristic polynomial of the Frobenius operation on the Monsky-Washnitzer cohomology, a well-behaved cohomology theory with
\mbox{$p$-adic}
coefficients for affine varieties in
characteristic~$p$.
This~approach seems to be very natural. It~is, however, still limited to curves, and in fact to particular types of~them. T.~Satoh's
\mbox{$p$-adic}
algorithm~\cite{Sa} for ordinary elliptic curves could perhaps be seen as a predecessor of this~method.

On~the other hand, as was shown by D.~Harvey~\cite{Ha}, an entirely elementary approach is possible, as well as a generalisation to arbitrary schemes of finite type over a finite~field. Harvey's method is, at least partially, based on earlier~ideas. We~are most likely unaware of some of its predecessors, but the work of D.~Wan~\cite{Wa99,Wa08}, as well as that of A.~Lauder and D.~Wan~\cite{LW}, certainly should be mentioned. Probably,~B.\ Dwork~\cite{Dw60a,Dw60b} is the very first, to whom the ideas behind~\cite{Ha} may be traced~back.

We~report on a variation of Harvey's
\mbox{$p$-adic}
method for double covers of the projective space, which we implemented as far as required for the needs of our search for
$K3$~surfaces
of degree two having RM or~CM. Although~it was never mentioned in print, this variation was certainly known, or at least obvious, to D.~Harvey~before. Generally~speaking, it seems that the algorithm described in~\cite{Ha}, although general in theory, may practically be implemented only in a form specialised to particular types of~varieties. Our code will be available within {\tt magma}~\cite{BCP}, version~2.22.\medskip

{\em Acknowledgement.}
We wish to thank David Harvey for valuable one-to-one talks, as well as the two anonymous referees for their suggestions on how to improve this~note.

\section{Counting points on varieties over finite fields--Elementary methods}

Let us assume that
$\smash{X \subset \Ab^N_{\bbF_{\!q}}}$
is an affine hypersurface. I.e., that
$X = V(f)$
for an arbitrary polynomial
$f \in \bbF_{\!q}[T_1,\ldots,T_N]$.
At~least in theory, this is not a serious restriction as every variety over a base field is birationally equivalent to a hypersurface~\cite[Chap.~I, Prop.~4.9]{Hart}, cf.\ \cite[\S1.3]{Ha}. In~practice, of course, too complicated transformations are undesirable, as they obviously reduce the speed of the~method. Even~if that concerns only a constant~factor.

\begin{ttt}[Determination of the points]
The most naive approach to point counting is certainly to {\em determine\/} all the points. If,~for whatever application, the actual points are asked for then there is no alternative to such a~method. Moreover,~to have some code available that realises a naive approach is useful for testing more advanced~algorithms.

In~order to determine the
$\bbF_{\!q^i}$-rational
points
on~$X$,
one has to run an iterated loop over all
$(N-1)$-tuples
$\smash{(x_1,\ldots,x_{N-1}) \in \bbF_{\!q^i}^{N-1}}$.
Each~time, the roots in
$\bbF_{\!q^i}$
of a univariate
polynomial~$g$
have to be found, which just means to compute
$\smash{\gcd(g, T^{q^i}-T)}$.
Thus,~the complexity to determine the
$\bbF_{\!q^i}$-rational
points is essentially
$O(q^{(N-1)i})$.
\end{ttt}

\begin{rem}
When~$X$
is defined over
$\bbF_{\!q}$,
but points over
$\bbF_{\!q^i}$
for
$i>1$
are sought for, one may gain a factor
of~$i$
by dealing with Frobenius orbits instead of~points.
\end{rem}

\begin{ttt}[FFT point counting]
Assume that
$X$
is given by a {\em decoupled\/}
polynomial~$f$.
I.e., that
$f = f_1(T_1,\ldots,T_M) + f_2(T_{M+1},\ldots,T_N)$
for some
$1 < M < N$.
Then~$\#X(\bbF_{\!q^i}) = (c_1 * c_2)(0)$
in terms of the two counting functions
$c_1$
and~$c_2$,
given by
\begin{align*}
c_1(c) & := \#\{(x_1,\ldots,x_M) \in \bbF_{\!q^i}^M \mid f_1(x_1,\ldots,x_M) = c \} \tag*{and\hspace{1.5cm}}\\
c_2(c) & := \#\{(x_{M+1},\ldots,x_N) \in \bbF_{\!q^i}^{N-M} \mid f_2(x_{M+1},\ldots,x_N) = c \} \, .
\end{align*}
In~order to compute
$\#X(\bbF_{\!q^i})$,
one first compiles look-up tables for
$c_1$
and~$c_2$
and then calculates
$(c_1 * c_2)(0)$,
according to the~definition. The complexity of this method to count the
$\bbF_{\!q^i}$-rational
points is
$O(q^{\max\{M,N-M\}i})$.

Similarly,~for the double cover of
$\Ab^N$
given by
$W^2 = f$,
one~has
$$\#X(\bbF_{\!q^i}) = q^{iN} + \sum_{\smash{c \in \bbF_{\!q^i}}} \!\chi(c) \!\cdot\! (c_1 * c_2)(c) \, ,$$
where~$\chi\colon \bbF_{\!q^i} \to \{-1,0,1\}$
denotes the quadratic~character. The~convolution
$c_1 * c_2$
then has to be computed using the FFT method~\cite[S\"atze~20.2 and~20.3]{Fo} or \cite[Theorem~32.8]{CLR}. Here,~again, the complexity is essentially
$O(q^{\max\{M,N-M\}i})$,
with an additional
\mbox{$\log$-factor}
coming from the Fourier~transforms.
Cf.~\cite[Algorithm~17]{EJants8} for more details on this~method.
\end{ttt}

\begin{ex}[cf.~{\cite[Example~3.13]{EJ15}}]
Consider the
$K3$~surface
$X$
of
degree~$2$
over~$\bbF_{\!7}$,
given~by
\begin{align*}
W^2 = 6 T_0^6 + 6T_0^5T_1 + 2T_0^5T_2 + 6T_0^4T_1^2 + 5T_0^4T_2^2 + 5T_0^3T_1^3 + T_0^2T_1^4 + 6T_0T_1^5 +{}& 5T_0T_2^5 \\[-1mm]
&{} + 3T_1^6 + 5T_2^6 \,.
\end{align*}
The polynomial on the right hand side does not contain any monomial involving both
$T_1$
and~$T_2$.
Thus,~over the affine plane
$\Pb^2_{\overline\bbF_{\!7}} \!\setminus\! V(T_0)$,
the double cover
$X$~is
given by a decoupled~polynomial.
The~numbers of points over
$\bbF_{7},\ldots,\bbF_{7^{10}}$
are
\begin{align*}
60,~
2\,488,~
118\,587,~
5\,765\,828,~
282\,498\,600,&~
13\,841\,656\,159,~
678\,225\,676\,496,\\
33\,232\,936\,342\,644,~
&1\,628\,413\,665\,268\,026,~
79\,792\,266\,679\,604\,918.
\end{align*}
Using~the FFT method, it took about 2 hours of CPU time and required approximately 5\,GB of memory to compute these numbers on an AMD 248 Opteron processor running at 2.2\,GHz.
\end{ex}

\begin{ttt}[Making \'etale cohomology explicit]
For a general variety, it is probably a tough ask to make its
\mbox{$l$-adic}
cohomology vector spaces explicit, including the Frobenius operation on~them. There~are, however, certain special types of varieties, for which this is possible. For~example, for an elliptic curve, one has
$\smash{H_\et^1(X_{\overline\bbF_{\!q}}, \bbQ_l) \cong \liminv_n E(\overline\bbF_{\!q})[l^n]} \!\otimes_{\bbZ_l}\! \bbQ_l$,
which is the starting point of Schoof's algorithm~\cite{Sch}.

A~different kind are those varieties, for which all the cycle maps~\cite[Cycle, \S2.2.10]{SGA4h}
$$\cl\colon \CH^j(X_{\overline\bbF_{\!q}}) \!\otimes_\bbZ\! \bbQ_l \longrightarrow H_\et^{2j}(X_{\overline\bbF_{\!q}}, \bbQ_l(j))$$
are~surjective and, moreover,
$H_\et^j(X_{\overline\bbF_{\!q}}, \bbQ_l) = 0$
for every odd
integer~$j$.
These~assumptions are fulfilled, for instance, for surfaces that are geometrically~rational.

A~case for which this approach is implemented in {\tt magma} is that of a cubic~surface.
Here,~$\CH^1(X_{\overline\bbF_{\!q}}) = \Pic(X_{\overline\bbF_{\!q}}) \cong \bbZ^7$
is generated by the 27~lines. One~has to determine the~lines by a Gr\"obner base calculation. From~the operation
of~$\Frob$
on the lines, one deduces the operation
on~$H_\et^2(X_{\overline\bbF_{\!q}}, \bbQ_l(1))$.
\end{ttt}

\begin{ex}
\leavevmode
\vspace{-.8mm}
{\footnotesize
\begin{verbatim}
> p := NextPrime(31^31);
> p;
17069174130723235958610643029059314756044734489
> rr<x,y,z,w> := PolynomialRing(GF(p),4);
> gl := x^3 + 2*y^3 + 3*z^3 + 5*w^3 - 7*(x+y+z+w)^3;
> time NumberOfPointsOnCubicSurface(gl);
291356705504951337929728147720395643095420441743546148074332970578622743757660955298550825611
Time: 1.180
\end{verbatim}}
\end{ex}

\begin{rem}
The~list of methods given in the section is by no means meant to be exhaustive. There~is at least one further idea, which we can mention only in~passing, namely the use of a~fibration. For~more information, we advise the reader to consult the articles \cite{La} of A.~Lauder and~\cite{PT} of S.~Pancratz and J.~Tuitman.
\end{rem}

\section{Harvey's
\mbox{$p$-adic}
method}

As~above, we assume that
$\bbF_{\!q}$
is a finite field of
characteristic~$p>0$.
The~method is
\mbox{$p$-adic}
in~nature. It~finds an approximation
of~$\#X(\bbF_{\!q^i})$
with respect to the
\mbox{$p$-adic}
valuation. Thus,~in~order to pinpoint
$\#X(\bbF_{\!q^i})$
exactly, one needs, in addition, an estimate for the number of~points.\medskip\smallskip

{\em Characteristic functions.}

\begin{nota}
\label{quasi}

\begin{iii}
\item
\looseness-1
We~continue to assume that
$X$
is an affine hypersurface. However,~we now suppose that
$\smash{X := V(f) \cap \bbG_{m,\bbF_{\!q}}^N \subset \bbG_{m,\bbF_{\!q}}^N}$
is a subscheme of an affine~torus. Again,~this restriction is not serious as
$\Ab^N_{\bbF_{\!q}}$
is a finite union of affine~tori. For~varieties more general than hypersurfaces, one might, in principle, simplify the situation by iteratively projecting away from points, cf.~\cite[\S1.3]{Ha}.
\item
For~every natural
number~$i$,
let us denote by
$\mu_{q^i-1}$
the group of all
\mbox{$(q^i-1)$-th}
roots~of unity and write
$\bbZ_{q^i}$
for the integer ring of the local field
$\bbQ_p(\mu_{q^i-1})$,
which is an unramified extension
of~$\bbQ_p$.
$\bbZ_{q^i}$~is
thus a complete discrete valuation ring with residue
field~$\bbF_{\!q^i}$.
Furthermore,~$\mu_{q^i-1} \subset \bbZ_{q^i}$
is bijectively mapped onto
$\smash{\bbF_{\!q^i}^*}$
under the residue map
$$\pi = \pi_i\colon \bbZ_{q^i} \twoheadrightarrow \bbF_{\!q^i} \, .$$
This~is, of course, an instance of taking Teichm\"uller representatives.
\item
Moreover,~one has
$\bbZ_{q^i} \subseteq \bbZ_{q^{i'}}$
if and only if
$\bbF_{\!q^i} \subseteq \bbF_{\!q^{i'}}$.
In~this case, we introduce the {\em quasi norm\/} map
\begin{eqnarray*}
N_{\bbZ_{q^{i'}}/\bbZ_{q^i}}\colon \bbZ_{q^i}[T_1,\ldots,T_N] &\longrightarrow& \bbZ_{q^i}[T_1,\ldots,T_N] \, , \\[-1mm]
g &\mapsto& g \!\cdot\! g^{(q^i)} \!\cdot\! g^{(q^{2i})} \cdots g^{(q^{i'-i})} \, ,
\end{eqnarray*}
for
$g^{(n)}$
the {\em quasi power,} given by
$g^{(n)}(T_1,\ldots,T_N) := g(T_1^n,\ldots,T_N^n)$.

The quasi norm is compatible with the usual norm map
$\smash{\N_{\bbF_{\!q^{i'}}/\bbF_{\!q^i}}\colon \bbF_{\!q^{i'}} \to \bbF_{\!q^i}}$
between the finite residue fields in the sense that, for arbitrary
$g \in \bbZ_{q^i}[T_1,\ldots,T_N]$
and
$z_1,\ldots,z_N \in \bbZ_{q^{i'}}$,
one~has
$N_{\bbZ_{q^{i'}}/\bbZ_{q^i}}g(z_1,\ldots,z_N) \in \bbZ_{q^i}$
and
$$\pi_i(N_{\bbZ_{q^{i'}}/\bbZ_{q^i}}g(z_1,\ldots,z_N)) = \N_{\bbF_{\!q^{i'}}/\bbF_{\!q^i}}(\overline{g}(\pi_{i'}(z_1),\ldots,\pi_{i'}(z_N))) \, .$$
Here,~$\overline{g} \in \bbF_{\!q^i}[T_1,\ldots,T_N]$
denotes the reduction
of~$g$
modulo~$(p)$.
\item
We~choose a lift
$F \in \bbZ_q[T_1,\ldots,T_N]$
of~$f$.
\end{iii}
\end{nota}

\begin{defi}
We~call a function
$\Phi\colon \mu_{q^i-1}^N \to \bbZ_q$
a {\em characteristic function\/}
of~$X(\bbF_{\!q^i})$
{\em modulo\/}~$p^l$
if, for every
$(z_1,\ldots,z_N) \in \mu_{q^i-1}^N$,
\begin{align*}
&\Phi(z_1,\ldots,z_N) \equiv 0 \pmod {p^l} \;\;\text{ if } (\pi(z_1),\ldots,\pi(z_N)) \not\in X(\bbF_{\!q^i}) \quad\quad\text{ and} \\
&\Phi(z_1,\ldots,z_N) \equiv 1 \pmod {p^l} \;\;\text{ if } (\pi(z_1),\ldots,\pi(z_N)) \in X(\bbF_{\!q^i}) \, .
\end{align*}
\end{defi}

\begin{prop}
\label{charfkt}
Let\/~$a$
be a positive integer such that\/
$a(q-1) \geq l$
and\/
$i \in \bbN$
be~arbitrary. Then
\begin{eqnarray*}
\Phi & := & (1-N_{\bbZ_{q^i}/\bbZ_q}F^{a(q-1)})^l \\
 & = & 1 + \sum\limits_{k=1}^l (-1)^k \binom{l}k N_{\bbZ_{q^i}/\bbZ_q}F^{ka(q-1)}
\end{eqnarray*}
is a characteristic function
of\/~$X(\bbF_{\!q^i})$
modulo\/~$p^l$.\medskip

\noindent
{\bf Proof.}
{\em
By~definition,
$(\pi(z_1),\ldots,\pi(z_N)) \in X(\bbF_{\!q^i})$
if and only if
$f(\pi(z_1),\ldots,\pi(z_N)) = 0 \in \bbF_{\!q^i}$,
which is equivalent to
$\N_{\bbF_{\!q^i}/\bbF_{\!q}}(f(\pi(z_1),\ldots,\pi(z_N))) = 0 \in \bbF_{\!q}$
and
$$N_{\bbZ_{q^i}/\bbZ_q}F(z_1,\ldots,z_N) \equiv 0 \pmod p \, .$$
If~this is true then
$\smash{N_{\bbZ_{q^i}/\bbZ_q}F^{a(q-1)}(z_1,\ldots,z_N) \equiv 0 \pmod {p^{a(q-1)}}}$,
hence the same modulo
$p^l$,
which implies
$\smash{1-N_{\bbZ_{q^i}/\bbZ_q}F^{a(q-1)}(z_1,\ldots,z_N) \equiv 1 \pmod {p^l}}$,
and
$$(1-N_{\bbZ_{q^i}/\bbZ_q}F^{a(q-1)}(z_1,\ldots,z_N))^l \equiv 1 \pmod {p^l} \, .$$
On~the other hand, in the opposite case, one has
$\smash{N_{\bbZ_{q^i}/\bbZ_q}F^{q-1}(z_1,\ldots,z_N) \equiv 1 \pmod p}$
and, consequently,
$1-N_{\bbZ_{q^i}/\bbZ_q}F^{a(q-1)}(z_1,\ldots,z_N) \equiv 0 \pmod p$,
which shows
$$(1-N_{\bbZ_{q^i}/\bbZ_q}F^{a(q-1)}(z_1,\ldots,z_N))^l \equiv 0 \pmod {p^l} \, , $$
as required.
\eop
}
\end{prop}

\begin{rems}

\begin{iii}
\item
According to the definition, if
$\Phi$
is a characteristic function
of~$X(\bbF_{\!q^i})$
modulo~$p^l$
then
\begin{equation}
\label{summ}
\#X(\bbF_{\!q^i}) \equiv \!\!\!\!\!\!\!\!\!\sum_{(z_1,\ldots,z_N) \in \mu_{q^i-1}^N}\!\!\!\!\!\!\!\!\!\!\! \Phi(z_1,\ldots,z_N) \pmod {p^l} \, . 
\end{equation}
\item
One~has
$$\sum_{z \in \mu_{q^i-1}} \!\!\!\!z^e = 0$$
unless
$e$
is a multiple
of~$q^i-1$.
Therefore,~the right hand side of~(\ref{summ}) is
$(q^i-1)^N$~times
the sum of the coefficients
of~$\Phi$
at the monomials of the form
$\smash{(T_1^{e_1} \cdots T_N^{e_N})^{q^i-1}}$,
for~$e_1, \ldots, e_N \in \bbZ_{\geq0}$.
\item
For~$l=1$,
these facts are known for a long time, cf.~\cite[Sec.~1.1, Theorem~4]{BS}.
\item
One might want to take the step from
$l=1$
to a larger value in a more naive manner as~follows.
For~$(\pi(z_1),\ldots,\pi(z_N)) \not\in X(\bbF_{\!q^i})$,
one has
$N_{\bbZ_{q^i}/\bbZ_q}F^{a(q-1)}(z_1,\ldots,z_N) = 1 + sp$,
for some
unknown~$s \in \bbZ_{q^i}$.
This~yields
\begin{align*}
N_{\bbZ_{q^i}/\bbZ_q}F^{2a(q-1)}(z_1,\ldots,z_N) & = 1 + 2sp +  s^2p^2 \, , \\
N_{\bbZ_{q^i}/\bbZ_q}F^{3a(q-1)}(z_1,\ldots,z_N) & = 1 + 3sp + 3s^2p^2 + s^3p^3 \, , \\[-2mm]
& \hspace{1.8mm}\vdots \\
N_{\bbZ_{q^i}/\bbZ_q}F^{la(q-1)}(z_1,\ldots,z_N) & = \textstyle 1 + lsp + \binom{l}2 s^2p^2 + \binom{l}3 s^3p^3 + \cdots + \binom{l}{l-1} s^{l-1}p^{l-1} + \binom{l}l s^lp^l \, .
\end{align*}
There~is a linear combination of these expressions that eliminates all terms containing
$p, p^2, \ldots, p^{l-1}$.
Namely,~one needs to take as the coefficient vector the first row of the Pascal-like matrix
$$
\left(
\begin{array}{ccccc}
1 & 1 & 0 & \cdots & 0 \\
1 & 2 & 1 & \cdots & 0 \\[-1mm]
\vdots & \vdots & \vdots & \ddots & \vdots \\
1 & l\!-\!1 & \binom{l-1}2 & \cdots & \binom{l-1}{l-1} \\
1 & l & \binom{l}2 & \cdots & \binom{l}{l-1}
\end{array}
\right)^{-1} .
$$
However,~well-known relations between binomial coefficients reveal that the first row of this matrix is exactly
$(l, -\binom{l}2, \binom{l}3, \ldots, (-1)^{l+1}\binom{l}l)$,
in agreement with Proposition~\ref{charfkt}.
\end{iii}
\end{rems}

{\em Linear operators.}

\begin{nota}

\begin{iii}
\item
We~consider the polynomial ring
$\bbZ_q[T_1,\ldots,T_N]$
as a free
\mbox{$\bbZ_q$-module}
and equip it with the symmetric bilinear form
$\langle.,.\rangle$,
given by
$$\left\langle \sum_{i_1,\ldots,i_N}\!\!\! a_{i_1,\ldots,i_N} T_1^{i_1}\cdots T_N^{i_N}, \sum_{i_1,\ldots,i_N}\!\!\! b_{i_1,\ldots,i_N} T_1^{i_1}\cdots T_N^{i_N} \right\rangle := \sum_{i_1,\ldots,i_N}\!\!\! a_{i_1,\ldots,i_N}b_{i_1,\ldots,i_N} \, .$$
Then~the monomials form an orthonormal~basis. Furthermore,~the bilinear form
$\langle .,. \rangle$
yields a
norm~$\|.\|$
on~$\bbZ_q[T_1,\ldots,T_N]$
such that
$\|g\| := [\|\langle g, g\rangle \|_p]^{1/2}$,
for
$\|.\|_p$
the usual normalised
\mbox{$p$-adic}~valuation.
The~completion
$\bbZ_q\{\{T_1,\ldots,T_N\}\}$
is a
\mbox{$\bbZ_q$-Hilbert}
space.
\item
There is the bounded linear operator
$\varphi\colon \bbZ_q\{\{T_1,\ldots,T_N\}\} \to \bbZ_q\{\{T_1,\ldots,T_N\}\}$,
provided by taking the quasi power
$h \mapsto h^{(q)}$.
Its~adjoint
$$\kappa\colon \bbZ_q\{\{T_1,\ldots,T_N\}\} \longrightarrow \bbZ_q\{\{T_1,\ldots,T_N\}\}$$
maps the monomial
$(T_1^{e_1}\cdots T_N^{e_N})^q$
to~$T_1^{e_1}\cdots T_N^{e_N}$
and all monomials not being of this particular type to~zero. In~particular,
$\kappa$
is a retraction for~$\varphi$.
I.e., one has
$\kappa \!\circ\! \varphi = \id$.

On~the other hand, for
$i \in \bbN$,
the composition
$\varphi^i \!\circ\! \kappa^i \colon \bbZ_q\{\{T_1,\ldots,T_N\}\} \to \bbZ_q\{\{T_1,\ldots,T_N\}\}$
is the orthogonal projection to the
sub-$\bbZ_q$-Hilbert
space spanned by all monomials of the type
$\smash{(T_1^{e_1}\cdots T_N^{e_N})^{q^i}}$.
\item
For~every polynomial
$g \in \bbZ_q[T_1,\ldots,T_N]$,
the multiplication map
$$m_g\colon \bbZ_q\{\{T_1,\ldots,T_N\}\} \longrightarrow \bbZ_q\{\{T_1,\ldots,T_N\}\} \, , \quad h \mapsto gh \, ,$$
is a bounded linear~operator.\vspace{-2mm}
\end{iii}
\end{nota}

\begin{fac}
\label{cmp}
Let\/~$g \in \bbZ_q[T_1,\ldots,T_N]$
and\/~$i \in \bbN$
be arbitrary. Then
$$\kappa^i \!\circ\! m_{g^{(q^i)}} \!\circ\! \varphi^i = m_g \, . \eop\vspace{-5mm}$$
\end{fac}

\begin{defi}[Linear operator associated with a polynomial]
\label{tr_op}
Given~a positive
integer~$i$
and a polynomial
$g \in \bbZ_q[T_1,\ldots,T_N]$,
we consider the linear operator
$$\M_{g,i} := \kappa^i \!\circ\! m_g \,.$$
\end{defi}

This operator is of interest because of the two results below.\vspace{-3mm}

\begin{lem}
For~every positive integer\/
$i$
and every
polynomial\/~$g \in \bbZ_q[T_0,\ldots,T_N]$,
one~has
\begin{equation}
\label{mat_prod}
\M_{\N_{\bbZ_{q^i}/\bbZ_q} \!g,\,i} = (\M_{g,1})^i \, .
\end{equation}
{\bf Proof.}
{\em
We~will show this inductively. The~assertion is clearly true for
$i=1$.
For~the inductive step, let us assume that equation~(\ref{mat_prod}) is true for a certain value
of~$i$.
Then,~by \ref{quasi}.iii), we have
$\N_{\bbZ_{q^{i+1}}/\bbZ_q} g = g^{(q^i)} \!\cdot\! \N_{\bbZ_{q^i}/\bbZ_q} g$.
Hence,~Definition~\ref{tr_op} shows~that
\begin{eqnarray*}
\M_{\N_{\bbZ_{q^{i+1}}/\bbZ_q} \!g,\,i+1} & = & \kappa^{i+1} \!\circ\! m_{\N_{\bbZ_{q^{i+1}}/\bbZ_q} \!g} \\
 & = & \kappa^{i+1} \!\circ\! m_{g^{(q^i)}} \!\circ\! m_{\N_{\bbZ_{q^i}/\bbZ_q} \!g} \\
 & = & \kappa^{i+1} \!\circ\! m_{g^{(q^i)}} \!\circ\! \varphi^i \!\circ\! \kappa^i \!\circ\! m_{\N_{\bbZ_{q^i}/\bbZ_q} \!g} \, .
\end{eqnarray*}
Observe~that the projection
$\varphi^i \!\circ\! \kappa^i$
has no effect here. Indeed,~it is followed by the multiplication with a polynomial, all of whose monomials have exponents only divisible
by~$q^i$,
and the
re\-trac\-tion~$\kappa^i$.
Thus,~according to Fact~\ref{cmp}, Definition~\ref{tr_op}, and the induction hypothesis,
\begin{eqnarray*}
\M_{\N_{\bbZ_{q^{i+1}}/\bbZ_q} \!g,\,i+1} & = & \kappa \!\circ\! m_g \!\circ\! \M_{\N_{\bbZ_{q^i}/\bbZ_q} \!g,\,i} \\
 & = & \M_{g,1} \!\circ\, (\M_{g,1})^i \\[1.4mm]
 & = & (\M_{g,1})^{i+1} \, ,
\end{eqnarray*}
as required.
\eop
}
\end{lem}

\begin{prop}
Let\/~$i$
be a positive integer and\/
$g \in \bbZ_q[T_0,\ldots,T_N]$
be any~polynomial.

\begin{abc}
\item
Then the linear operator\/
$\M_{g,i}$
is trace~class.
\item
Its trace\/
$\tr \M_{g,i}$
is equal to the sum of all coefficients
of\/~$g$
at the monomials of the form\/
$(T_1^{e_1} \cdots T_N^{e_N})^{q^i-1}$.
\end{abc}\smallskip

\noindent
{\bf Proof.}
{\em
We~will only prove b), as this immediately implies~a). By~definition,
$$\tr \M_{g,i} = \sum_{k\in K} \langle \M_{g,i}(v_k), v_k \rangle \, ,$$
for
$(v_k)_{k\in K}$
any complete orthonormal~system. In~particular, for the system of all~monomials.

For~the monomial
$T_1^{e_1} \cdots T_N^{e_N}$,
the pairing
$\langle \M_{g,i}(T_1^{e_1} \cdots T_N^{e_N}), T_1^{e_1} \cdots T_N^{e_N} \rangle$
yields the coefficient of
$\M_{g,i}(T_1^{e_1} \cdots T_N^{e_N}) = \kappa^i (gT_1^{e_1} \cdots T_N^{e_N})$
at the monomial
$T_1^{e_1} \cdots T_N^{e_N}$.
This~is nothing but the coefficient of
$gT_1^{e_1} \cdots T_N^{e_N}$
at
$\smash{(T_1^{e_1} \cdots T_N^{e_N})^{q^i}}$
and, hence, the coefficient
of~$g$
at
$\smash{(T_1^{e_1} \cdots T_N^{e_N})^{q^i-1}}$.
The assertion~follows.
}
\eop
\end{prop}

\begin{rem}
For~general information on
$p$-adic
Hilbert spaces, linear operators, and the concept of the trace, the interested reader might compare Chapters 2 and~3 of the textbook~\cite{Di} of T.~Diagana.
\end{rem}

{\em An algorithm.}

\begin{algo}
\label{Harvey_alg}
Given a projective hypersurface
$X = V(f) \subset \Pb^N_{\bbF_{\!q}}$
and a positive
integer~$n$,
this algorithm computes the numbers of points
$\#X(\bbF_{\!q})$,
$\#X(\bbF_{\!q^2})$,
\ldots,~$\#X(\bbF_{\!q^n})$.

\begin{iii}
\item
Choose a lift
$\smash{\widetilde{F} \in \bbZ_q[T_0,\ldots,T_N]}$
of~$f$.
Decompose~$X$
into~$X^c := [V(T_0) \cup\ldots\cup V(T_N)] \cap X$
and~$X^0 := X \!\setminus\! X^c$.
For~the lower dimensional closed subscheme
$X^c$,
either use naive point counting methods or apply the method recursively.
For~$X^0$,
do the following.
\item
Use elementary arguments or the Weil conjectures or both to find estimates for
$\#X^0(\bbF_{\!q})$,
$\#X^0(\bbF_{\!q^2})$,
\ldots,~$\#X^0(\bbF_{\!q^n})$.
From~these, deduce a
\mbox{$p$-adic}
precision~$l$
that suffices to determine these numbers~exactly. Put,~finally,
$a := \lceil l/(q-1) \rceil$.
\item
Define
$F \in \bbZ_q[T_1,\ldots,T_N]$
by putting
$\smash{F(T_1,\ldots,T_N) := \widetilde{F}(1,T_1,\ldots,T_N)}$.
\item
Find matrix representations of the operators
$\M_{F^{a(q-1)},1}$,
$\M_{F^{2a(q-1)},1}$,
\ldots,
$\M_{F^{la(q-1)},1}$.
Work~on finite dimensional subspaces that are large enough to get the traces~right. For~a polynomial of total
degree~$d$,
as a basis one may use all monomials of total
degree~$\leq\!\! \lceil d/(q-1) \rceil$.
\item
For~$i := 1$
to~$n$,
do the following.

Calculate the matrix powers
$\M_{F^{a(q-1)},1}^i$,
$\M_{F^{2a(q-1)},1}^i$,
\ldots,
$\M_{F^{la(q-1)},1}^i$
and eventually their traces
$\smash{\tr \M_{F^{ka(q-1)},1}^i}$.
Then
\begin{equation}
\label{point_count}
\#X^0(\bbF_{\!q^i}) \equiv (q^i-1)^N \left[ 1 + \sum\limits_{k=1}^l (-1)^k \binom{l}k \tr \M_{F^{ka(q-1)},1}^i \right] \pmod {p^l} \, .
\end{equation}
Determine the exact value of
$\#X^0(\bbF_{\!q^i})$
from this congruence and the estimate established in step~ii).
Add~$\#X^c(\bbF_{\!q^i})$
and output the integer~found.
\end{iii}
\end{algo}

\begin{rems}

\begin{iii}
\item
Formula~(\ref{point_count}) is essentially the same as the ``trace formula'' in \cite[Theorem~3.1]{Ha}.
\item[ii) (Cohomological interpretation) ]
The matrix
$\M_{F^{q-1},1}$
appears for the first time in \cite[\S2]{Dw60a} and later in \cite[Corol\-lar\-ies~1 and~2]{Mi}, where it is shown to be the natural generalisation of the classical Hasse-Witt matrix~\cite{HW}. In~other words, it represents the Frobenius homomorphism on the middle coherent cohomology
$H^{N-1}(X,\calO_X)$.
The~characteristic polynomial
of~$\M_{F^{q-1},1}$
is hence congruent 
modulo~$p$
to the Weil polynomial
$\chi_{N-1}$
\cite[Exp.\,XXII, Th\'eor\`eme~3.1]{SGA7}. Unfortunately,~the algorithm above does not provide a canonical lift representing the crystalline Frobenius, but only successive
\mbox{$p$-adic}
approximations of the proper eigenvalues, so that (\ref{point_count}) finally holds due to many~cancellations.
\end{iii}
\end{rems}

\section{A variation adapted to double covers}

{\em The variation.}

\begin{nota}

\begin{iii}
\item
We~assume in this section that
$\bbF_{\!q}$
is a finite field of
characteristic~$p \neq 2$.
\item
Furthermore,~let
$X$
be a double cover of an affine torus
$\smash{\bbG_{m,\bbF_{\!q}}^N}$,
given by an equation of the form
$w^2 = f$
for a polynomial
$f \in \bbF_{\!q}[T_1,\ldots,T_N]$.
Our~interest stems from double covers of projective space
of~$\smash{\Pb_{\bbF_{\!q}}^N}$,
but, as before, going over from projective space to an affine torus is not a serious~modification. One~may at the very end cover
$\smash{\Pb_{\bbF_{\!q}}^N}$
by affine~tori.
\item
We~choose a lift
$F \in \bbZ_q[T_1,\ldots,T_N]$
of~$f$.
\item
For a positive
integer~$l$
and
$k = 1,\ldots,l$,
let us put
\begin{equation}
\label{ali}
A^{(l)}_k := \frac{(2l)!}{2^{2l-1}l!(l-1)!} \!\cdot\! \frac{(-1)^{k+1}}{2k-1} \binom{l-1}{k-1} \, .
\end{equation}
\end{iii}
\end{nota}

\begin{lem}
\label{Pascal}
Let\/~$l > 0$
and\/~$b$
be~integers.

\begin{abc}
\item
Then~the Pascal-like~matrix
\begin{equation}
\label{Pasc_mat}
P_b^{(l)} :=
\left(
\begin{array}{lccccc}
1 & b & \binom{b}2 & \binom{b}3 &\cdots & \binom{b}{l-1} \\
1 & b+2 & \binom{b+2}2 & \binom{b+2}3 &\cdots & \binom{b+2}{l-1} \\
\vdots & \vdots & \vdots & \vdots & \ddots & \vdots \\
1 & b\!+\!2l\!-\!4 & \binom{b\!+\!2l\!-\!4}2 & \binom{b\!+\!2l\!-\!4}3 &\cdots & \binom{b\!+\!2l\!-\!4}{l-1} \\
1 & b\!+\!2l\!-\!2 & \binom{b\!+\!2l\!-\!2}2 & \binom{b\!+\!2l\!-\!2}3 &\cdots & \binom{b\!+\!2l\!-\!2}{l-1}
\end{array}
\right) \in M_{l\times l}(\bbQ)
\end{equation}
is~invertible.
\item
All coefficients of its inverse\/
$\smash{\big( P_b^{(l)} \big)^{-1}}$
are\/
\mbox{$p$-adic}
integers.
\item
The first row of\/
$\smash{\big( P_1^{(l)} \big)^{-1}}$
is\/
$(A^{(l)}_1, A^{(l)}_2, \ldots, A^{(l)}_l)$.
\end{abc}\smallskip

\noindent
{\bf Proof.}
{\em
a) and~b)
Column transformations convert
$\smash{P_b^{(l)}}$
essentially into a Vandermonde matrix. This~shows
$\smash{\det P_b^{(l)} = \frac{(2l-2)!!(2l-4)!!\ldots 2!!}{1!2!\ldots (l-1)!} = 2^{(l-1)+(l-2)+\cdots+1} = 2^{\binom{l}2} \neq 0}$,
which immediately implies~a). Assertion~b) follows from this together with Cramer's~rule.\smallskip

\noindent
c)
The assertion is that the linear combination with coefficients
$(A^{(l)}_1, A^{(l)}_2, \ldots, A^{(l)}_l)$
of the rows of the matrix
$\smash{P_1^{(l)}}$
yields~$(1, 0, \ldots, 0)$.
Here,~the occurrence of the zeroes is a consequence of the standard relation
$\sum_{k=1}^l (-1)^k \binom{l-1}{k-1} k^e = 0$
among binomial coefficients,
for~$e = 0, \ldots, l-2$.
It~is a trivial, but tedious work to adjust the constant factor.
}
\eop
\end{lem}

\begin{defi}
We~call a function
$\Phi\colon \mu_{q^i-1}^N \to \bbZ_q$
an {\em Eulerian function\/}
for~$X(\bbF_{\!q^i})$
{\em modulo\/}~$p^l$
if, for every
$(z_1,\ldots,z_N) \in \mu_{q^i-1}^N$,
\begin{align*}
&\Phi(z_1,\ldots,z_N) \equiv \phantom{-}1 \pmod {p^l} \quad\text{ if } f(\pi(z_1),\ldots,\pi(z_N)) \text{ is a square in } \bbF_{\!q^i}^* \, , \\
&\Phi(z_1,\ldots,z_N) \equiv -1 \pmod {p^l} \quad\text{ if } f(\pi(z_1),\ldots,\pi(z_N)) \text{ is a non-square in } \bbF_{\!q^i}^* \, , \text{ and}\\
&\Phi(z_1,\ldots,z_N) \equiv \phantom{-}0 \pmod {p^l} \quad\text{ if } f(\pi(z_1),\ldots,\pi(z_N)) = 0 \, .
\end{align*}
\end{defi}

\begin{theo}
Suppose that\/
$p$
is an odd prime and that\/
$b$
is odd such that\/
$b(q-1)/2 \geq l$.
Put\/~$A^{(l)}_{b,k} := \big(\big( P_b^{(l)} \big)^{-1}\big)_{1,k}$.
Then
$$\Phi := \sum\limits_{k=1}^l A^{(l)}_{b,k} \!\cdot\! N_{\bbZ_{q^i}/\bbZ_q}F^{(b+2k-2)(q-1)/2}$$
is an Eulerian function
for\/~$X(\bbF_{\!q^i})$
modulo\/~$p^l$.\medskip

\noindent
{\bf Proof.}
{\em
Recall~that each of the coefficients
$\smash{A^{(l)}_{b,1} := \big(\big( P_b^{(l)} \big)^{-1}\big)_{1,1}}$,
\ldots,
$\smash{A^{(l)}_{b,l} := \big(\big( P_b^{(l)} \big)^{-1}\big)_{1,l}}$
is a
\mbox{$p$-adic}
integer. Furthermore, one has that
$f(\pi(z_1),\ldots,\pi(z_N)) = 0 \in \bbF_{\!q^i}$
if and only if
$\smash{\N_{\bbF_{\!q^i}/\bbF_{\!q}}(f(\pi(z_1),\ldots,\pi(z_N))) = 0 \in \bbF_{\!q}}$,
which, in turn, is equivalent to
$\smash{N_{\bbZ_{q^i}/\bbZ_q}F(z_1, \ldots, z_N)}$
being a non~unit
in~$\bbZ_q$.
If~this is true then
$\smash{N_{\bbZ_{q^i}/\bbZ_q}F^{(b+2k-2)(q-1)/2}(z_1,\ldots,z_N)}$
is of
\mbox{$p$-adic}
valuation at least
$(b+2k-2)(q-1)/2 \geq b(q-1)/2 \geq l$.
Consequently,
$$\sum_{k=1}^l A^{(l)}_{b,k} \!\cdot\! N_{\bbZ_{q^i}/\bbZ_q}F^{(b+2k-2)(q-1)/2} \equiv 0 \pmod {p^l} \, .$$

Otherwise,
$\smash{f(\pi(z_1),\ldots,\pi(z_N)) \in (\bbF_{\!q^i}^*)^2}$
precisely when
$\smash{\N_{\bbF_{\!q^i}/\bbF_{\!q}}(f(\pi(z_1),\ldots,\pi(z_N))) \in \bbF_{\!q}^*}$
is a square. Then
$$\N_{\bbF_{\!q^i}/\bbF_{\!q}}(f^{(q-1)/2}(\pi(z_1),\ldots,\pi(z_N))) = 1 \in \bbF_{\!q} \, ,$$
while the same expression is equal
to~$(-1)$
in the case of a non-square. By~construction, this~means
$\smash{N_{\bbZ_{q^i}/\bbZ_q}F^{(q-1)/2}(z_1, \ldots, z_N) \equiv \pm1 \pmod p}$,
respectively.

Write
$N_{\bbZ_{q^i}/\bbZ_q}F^{(q-1)/2}(z_1, \ldots, z_N) = \pm(1 + sp)$
for some
unknown~$s \in \bbZ_{q^i}$.
This~yields
\begin{eqnarray*}
N_{\bbZ_{q^i}/\bbZ_q}F^{3(q-1)/2}(z_1,\ldots,z_N) & = & \pm\big(1 + 3sp + 3s^2p^2 + s^3p^3\big) \, ,\\
N_{\bbZ_{q^i}/\bbZ_q}F^{5(q-1)/2}(z_1,\ldots,z_N) & = & \pm\big(1 + 5sp + 10s^2p^2 + 10s^3p^3 + 5s^4p^4 + s^5p^5\big) \, , \\[-2mm]
&\vdots& \\
N_{\bbZ_{q^i}/\bbZ_q}F^{(b+2l-2)(q-1)/2}(z_1,\ldots,z_N) & = & \textstyle \pm\big(1 + (b\!+\!2l\!-\!2)sp + \binom{b+2l-2}2 s^2p^2 + \\[-1mm]
 & & \textstyle \hspace{1.45cm} \binom{b+2l-2}3 s^3p^3 + \cdots + \binom{b+2l-2}{l-1} s^{l-1}p^{l-1} + \ldots \big) \, .
\end{eqnarray*}
Lemma~\ref{Pascal}.a) and~b) now show
$\smash{\sum\limits_{k=1}^l A^{(l)}_{b,k} \!\cdot\! N_{\bbZ_{q^i}/\bbZ_q}F^{(b+2k-2)(q-1)/2}(z_1,\ldots,z_N) \equiv \pm1 \pmod {p^l}}$,
as required.%
}%
\eop
\end{theo}

\begin{rem}
According to the definition, if
$\Phi$
is an Eulerian function
for~$X(\bbF_{\!q^i})$
modulo~$p^l$
then
\begin{equation}
\#X(\bbF_{\!q^i}) \equiv (q^i-1)^N + \!\!\!\!\!\!\!\!\!\!\!\sum_{(z_1,\ldots,z_N) \in \mu_{q^i-1}^N}\!\!\!\!\!\!\!\!\!\!\! \Phi(z_1,\ldots,z_N) \pmod {p^l} \, .
\end{equation}
This leads to the algorithm~below.
\end{rem}

\begin{algo}
Given a polynomial
$f \in \bbF_{\!q}[T_0,\ldots,T_N]$
and a positive
integer~$n$,
this algorithm computes the numbers of points
$\#X(\bbF_{\!q})$,
$\#X(\bbF_{\!q^2})$,
\ldots,~$\#X(\bbF_{\!q^n})$
on the double cover
of~$\smash{\Pb^N_{\bbF_{\!q}}}$,
given
by~$w^2 = f$.

\begin{iii}
\item
Choose a lift
$\smash{\widetilde{F} \in \bbZ_q[T_0,\ldots,T_N]}$
of~$f$.
Decompose~$\smash{\Pb^N_{\bbF_{\!q}}}$
into~$P^c := V(T_0) \cup\ldots\cup V(T_N)$
and~$\smash{P^0 := \Pb^N_{\bbF_{\!q}} \!\setminus\! P^c}$.
Then~$\smash{P^0 \cong \bbG_{m,\bbF_{\!q}}^N}$.
For~the
restriction~$X^c$
of the double cover to the lower dimensional closed
subscheme~$P^c$,
either use naive point counting methods or apply the method recursively.
For~$X^0$,
the restriction
to~$P^0$,
do the~following.
\item
Use elementary arguments or the Weil conjectures or both to find estimates for
$\#X^0(\bbF_{\!q})$,
$\#X^0(\bbF_{\!q^2})$,
\ldots,~$\#X^0(\bbF_{\!q^n})$.
From~these, deduce a
\mbox{$p$-adic}
precision~$l$
that suffices to determine these numbers~exactly. Put, finally,
$b := \lfloor 2l/(q-1)\rfloor$.
\item
If~$b > 1$
then determine the coefficient vector
$\smash{(A^{(l)}_{b,1}, \ldots, A^{(l)}_{b,l})}$
by inverting the matrix
$P_b^{(l)}$
from~(\ref{Pasc_mat}). Otherwise,~calculate first
$$C_l := \frac{(2l)!}{2^{2l-1}l!(l-1)!}$$
\mbox{and then the coefficients
$\smash{\!A^{(l)}_k \!:=\! C_l \!\cdot\! \frac{(-1)^{k+1}}{2k-1} \binom{l-1}{k-1}}$.
Finally,
set~$\smash{\!(A^{(l)}_{1,1}, \ldots, A^{(l)}_{1,l}) \!:=\! (A^{(l)}_1, \ldots, A^{(l)}_l)}$.}
\item
Define
$F \in \bbZ_q[T_1,\ldots,T_N]$
by putting
$\smash{F(T_1,\ldots,T_N) := \widetilde{F}(1,T_1,\ldots,T_N)}$.
\item
\!Find \!matrix \!representations \!of \!the \!operators
$\M_{F^{(q-1)/2},1}$,
$\M_{F^{3(q-1/2)},1}$,
\ldots,
$\M_{F^{(b+2l-2)(q-1/2)},1}$.
Work~on finite dimensional subspaces that are large enough to get the traces~right. For~a polynomial of total
degree~$d$,
as a basis one may use all monomials of total
degree~$\leq\!\! \lceil d/(q-1) \rceil$.
\item
For~$i := 1$
to~$n$,
do the following.

Calculate the matrix powers
$\M_{F^{(q-1)/2},1}^i$,
$\M_{F^{3(q-1)/2},1}^i$,
\ldots,
$\M_{F^{(b+2l-2)(q-1)/2},1}^i$
and eventually their traces
$\smash{\tr \M_{F^{(b+2k-2)(q-1)/2},1}^i}$.
Then
$$\#X^0(\bbF_{\!q^i}) \equiv (q^i-1)^N \left[ 1 + \sum\limits_{k=1}^l A^{(l)}_{b,k} \!\cdot\! \tr \M_{F^{(b+2k-2)(q-1)/2},1}^i \right] \pmod {p^l} \, .$$
Determine the exact value of
$\#X^0(\bbF_{\!q^i})$
from this congruence and the estimate established in step~ii).
Add~$\#X^c(\bbF_{\!q^i})$
and output the integer~found.
\end{iii}
\end{algo}

\begin{rems}

\begin{iii}
\item (Complexity.)
The algorithm contains two time consuming operations. The~first is the computation of the powers
of~$F$
and the second is the multiplication of~matrices. We~tested our implementation in {\tt magma} on one core of an Intel i5-4690 processor running at~3.5\,GHz.

For~a
$K3$~surface
of degree two, in order to compute the whole Weil polynomial, the algorithm has to work with a
\mbox{$p$-adic}
precision of at least 11~digits. For~this, it has to compute
$F^{k(p-1)/2}$
for all odd integers
$1 \leq k \leq 21$.
This~results in matrices of size up to
$2080 \times 2080$
having entries
in~$\bbZ/p^{11}\bbZ$.

The table below lists the time used for a randomly chosen surface in a test modulo four representative primes of good~reduction. The computation required about 13\,GB of memory.
%
\begin{table}[H]
\small
\centerline{
\begin{tabular}{c|c|c|c}
$p$ & powers of $F$ & matrix build & matrix operations \\
\hline
 31  &  12.65 & 21.05 & 55.36 \\
 61  &  76.91 & 21.63 & 71.00 \\
 97  & 236.92 & 22.30 & 73.53 \\
127  & 489.92 & 22.36 & 73.97
\end{tabular}}
\caption{Time (in s) for the steps of the algorithm}
\end{table}\vspace{-3mm}

The~figures show that the time for the computation of the powers
of~$F$
is  approximately cubic
in~$p$.
This~is explained by the fact that the number of terms of the final result is quadratic
in~$p$,
while the number of multiplications necessary to approach it is linear
in~$p$.
The time for the matrix operations is almost constant as the size of the
matrices is constant. The variation of the time is explained completely by the
change of the ring the matrix entries are taken~from.

In~theory, there exists an algorithm that is quadratic
in~$p$.
This~is, however, not yet implemented in {\tt magma} and, reportedly, it is not easy to~implement.
\item
We~incorporated the variable
$b$
only for completeness. In~practice, the algorithm is of interest particularly when
$q$
is not very~small. Then
$q-1 \geq 2l$
and~$b = 1$.
The~same applies to the variable
$a$
in Algorithm~\ref{Harvey_alg}.
\item
One might as well determine the coefficient vector 
$\smash{(A^{(l)}_{b,1}, \ldots, A^{(l)}_{b,l})}$
by a matrix inversion, also when
$b=1$.
This~step is in fact not time-critical. Nevertheless, we find it interesting that it is possible to give an explicit~formula.
\end{iii}
\end{rems}

\begin{rem}[The moving simplex idea]
The~idea behind most of the operator calculus above is that only very few of the coefficients of the powers
of~$F$
are actually~needed. An~extreme case is the~following. Assume~that a
$K3$~surface
$X$
is given as a double cover
of~$\Pb^2_\bbQ$,
by the equation
$W^2 = F(T_0,T_1,T_2)$
for~$F$
a homogeneous polynomial of
degree~$6$.
Moreover,~one only wants to decide for which primes
$p$
up to a certain
bound~$B$
the reduction
$X_p$
is ordinary. This~means that only
$(\#X_p(\bbF_{\!p}) \bmod p)$
is asked~for.
Thus,~for each prime, just a single coefficient is to be computed, that of
$\smash{F^{\frac{p-1}2}}$
at
$T_0^{p-1} T_1^{p-1} T_2^{p-1}$.

\looseness-1
This~might be done as~follows. Consider~just a small triangle of coefficients
of~$F^e$,
those at
$T_0^{a_0} T_1^{a_1} T_2^{a_2}$
such that
$\smash{|a_i-\frac{a_0+a_1+a_2}3| \leq c_i}$
for
$i = 0, 1, 2$
and some constants
$c_0$,
$c_1$,
and~$c_2$.
These triangles may be calculated inductively for
$e = 1, \ldots, (B-1)/2$.
For~$F$
homogeneous of
degree~$d$,
any choice such that
$c_0+c_1+c_2 := 2(d-1)$
should~work.

\looseness-1
Indeed,~$F$
is given and therefore
$\frac{dF}{dT_0}$,
$\frac{dF}{dT_1}$,
and~$\frac{dF}{dT_2}$
are~known. When~the triangle for
$F^e$
has been established, computing
$F^{e+1}$
naively one would lose some of the~coefficients. But~the coefficients of
$F^{e+1}$
have massive linear~relations. A~single coefficient of
$F^{e+1}$
appears as a linear combination of the coefficients
of~$F^e$,
when evaluating the formula
$F^{e+1} = F \!\cdot\! F^e$.
It~appears as well in
$\smash{\frac{dF^{e+1}}{dT_0} \!=\! (e+1) F^e \frac{dF}{dT_0}}$
and the analogous formulae involving
$T_1$
and~$T_2$.
Hence,~each hypothetical coefficient yields three linear~relations. One~may thus reconstruct some neighbouring coefficients using the linear relations and therefore keep the size of the triangle constant under~iteration. This~approach stems from \cite[\S4.1]{Ha} and is implemented in the function {\tt NonOrdinaryPrimes}.

\looseness-1
On~the other hand, {\tt WeilPolynomialOfDegree2K3Surface} does not use the moving simplex idea, but the approach described above in this~section. The~point is that, for example for
$p = 127$,
one has to compute
$F^{63 \cdot 21}$,
which results in a sum of 31\,517\,830~terms. Thereby,~from the coefficients, about 14\% are used as matrix~entries. Thus, computing only the used coefficients of the powers
of~$F$
would speed up the algorithm only for considerably larger values
of~$p$.\vspace{-3mm}
\end{rem}

\begin{rem}
There are improvements to the moving simplex idea, which we did not implement. The point is that instead of keeping track of the coefficients
of~$F^e$
on the entire simplex, one may sometimes, and probably often, get by with working only on a relatively small~subset. For~example, the idea has been worked out for smooth space quartics, in which case 64 coefficients suffice out of the 220, one would naively~have. This~is called the controlled reduction Abbott--Kedlaya--Roe algorithm \cite[Algorithm~3.4.10]{AKR}. Cf.~\cite{Co}.
\end{rem}

\section{Families of
$K3$~surfaces
that are highly likely to have real or complex multiplication}
\label{RMCM}

\begin{ttt}
\label{arith_cons}
Based~on the Mumford-Tate conjecture, which was proven for
$K3$~surfaces
by S.\,G.~Tankeev~\cite{Ta90,Ta95}, we established in~\cite{EJ14} certain arithmetic consequences of real or complex multiplication for
$K3$~surfaces
defined
over~$\bbQ$.
Assume~for simplicity that the endomorphism field
$E$
is a quadratic number~field. Then~these consequences include the~following.

\begin{iii}
\item
If~$p$
is inert
in~$E$
then
$\#X_p(\bbF_{\!p}) \equiv 1 \pmod p$.
I.e.~$X_p$
is non-ordinary.
\item
For every prime
$p$
of good reduction, the transcendental factor
$\chi_p^\tri$
of the Weil polynomial
of~$X_p$
either splits over
$E$
into two factors conjugate to each other or becomes a square under raising all its roots to the
\mbox{$f$-th}
power, for
some~$f > 0$.
\end{iii}
During the last few years, we systematically searched for
$K3$~surfaces,
the reductions of which show such an unusually regular behaviour. So~far, we have the following conjectural list of suspicious~surfaces.
\end{ttt}

\begin{cons}
\label{list_ex}
Consider~the following\/
$K3$~surfaces
and families of~such,
\begin{align*}
V_a^{(2)} \colon
    & W^2 = \textstyle
 [(\frac18 a^2 \!-\! \frac12 a \!+\! \frac14)T_1^2 + (a^2 \!-\! 2a \!+\! 2)T_1T_2 + (a^2 \!-\! 4a \!+\! 2)T_2^2] \\[-.5mm]
    &  \textstyle \hspace{1.7mm}
 [(\frac18 a^2 \!+\! \frac12 a \!+\! \frac14)T_0^2 + (a^2 \!+\! 2a \!+\! 2)T_0T_2 + (a^2 \!+\! 4a \!+\! 2)T_2^2]
 [2T_0^2 + (a^2 \!+\! 2)T_0T_1 + a^2T_1^2] , \\
V_a^{(5)} \colon
    & W^2 = \textstyle
 [T_1^2 + aT_1T_2 + (\frac5{16} a^2 \!+\! \frac54 a \!+\! \frac54)T_2^2]
 [T_0^2 + T_0T_2 + (\frac1{320} a^2 \!+\! \frac1{16} a \!+\! \frac5{16}) T_2^2] \\[-.5mm]
    &   \textstyle\hspace{97.2mm}
 [T_0^2 + T_0T_1 + \frac1{20} T_1^2] , \\
V^{(13)} \colon
    & W^2 = (25T_1^2 + 26T_1T_2 + 13T_2^2)(T_0^2 + 2T_0T_2 + 13T_2^2)(9T_0^2 + 26T_0T_1 + 13T_1^2) \, , \\
V_{a,b}^{(2)} \colon & W^2 = T_0 T_1 T_2 f_{a,b} \quad\text{\rm for} \\[-.5mm]
    &   \hspace{1cm} f_{a,b} \!:= a (-T_0^2 T_1 \!+\! T_0^2 T_2 \!+\! 2 T_0 T_1^2 \!-\! 3 T_0 T_1 T_2 \!+\! T_0 T_2^2 \!+\! T_1^3 \!-\! 4 T_1^2 T_2 \!+\! 5 T_1 T_2^2 \!-\! 2 T_2^3) \\[-.5mm]
    &   \hspace{2cm} {}+ b (T_0^3 - 2 T_0^2 T_1 - T_0 T_1^2 + 3 T_0 T_1 T_2 - T_0 T_2^2 + T_1^2 T_2 - T_1 T_2^2) \, , \\
V_{a,b}^{(3)} \colon & W^2 = T_0 T_1 T_2 f_{a,b} \quad\text{\rm for} \\[-.5mm]
    &   \hspace{1cm} f_{a,b} := (T_0^3 - 2 T_0^2 T_2 - T_0 T_1^2 - T_0 T_2^2 - 2 T_1^2 T_2 + 2 T_2^3) a^2  \\[-.5mm]
    &   \hspace{2cm} {}+ (6 T_0^2 T_1 + 6 T_0^2 T_2 + 6 T_0 T_1^2 + 6 T_1^2 T_2 - 6 T_2^3) a b \\[-.5mm]
    &   \hspace{2cm} {}+ (-3 T_0^2 T_1 - 6 T_0^2 T_2 + 3 T_1^3 - 6 T_1^2 T_2 - 3 T_1 T_2^2 + 6 T_2^3) b^2 , \\
V_{{\bf a},{\bf b}}^{(-1)} \colon
    & W^2 = T_0 T_1 T_2 (T_0+T_1+T_2)(a_1 T_0 + a_2 T_1 + a_3 T_2)(b_1 T_0 + b_2 T_1 + b_3 T_2) \\[-.8mm]
    &   \hspace{1.3cm} \text{\rm for } {\bf a}, {\bf b} \in \smash{\bbC^3} \text{\rm ~such that } a_1 b_3 \!+\! a_2 b_1 \!-\! 2 a_3 b_1 = 0 \text{\rm ~and } a_1 b_2 \!+\! a_2 b_3 \!-\! 2 a_3 b_2 = 0 \, . \\
V^{(-1,\mu_7)} \colon
    & W^2 = T_0T_1T_2(7 T_0^3 - 7 T_0^2 T_1 + 49 T_0^2 T_2 - 21 T_0 T_1 T_2 + 98 T_0 T_2^2 + T_1^3 - 7 T_1^2 T_2 \!+\! 49 T_2^3) \, , \\
V^{(-1,\mu_9)} \colon
    & W^2 = T_0 T_1 T_2 (T_0^3 - 3 T_0^2 T_2 - 3 T_0 T_1^2 - 3 T_0 T_1 T_2 + T_1^3 + 9 T_1^2 T_2 + 6 T_1 T_2^2 + T_2^3) \, , \\
V^{\!(-1,\mu_{19})} \colon
    & W^2 = T_0 T_1 T_2 (49T_0^3 - 304T_0^2T_1 + 570T_0^2T_2 + 361T_0T_1^2 - 2793T_0T_1T_2 + 2033T_0T_2^2 \\[-.5mm]
    & \hspace{57mm} {} + 361T_1^3 + 2888T_1^2T_2 - 5415T_1T_2^2 + 2299T_2^3) \, .
\end{align*}%
\begin{abc}
\item
Then the generic fibre of each family, as well as each of the four individual surfaces has geometric Picard
rank\/~$16$.
\item
Moreover,~the generic fibre
of\/
$\smash{V_a^{(2)}}$
has real multiplication
by\/
$\bbQ(\sqrt{2})$,
that
of\/
$\smash{V_a^{(5)}}$
RM
by\/
$\bbQ(\sqrt{5})$,
that
of\/
$\smash{V_{a,b}^{(2)}}$
RM
by\/
$\bbQ(\sqrt{2})$,
that
of\/
$\smash{V_{a,b}^{(3)}}$
RM
by\/
$\bbQ(\sqrt{3})$,
and that
of\/
$\smash{V_{{\bf a},{\bf b}}^{(-1)}}$
complex multiplication
by\/~$\bbQ(\sqrt{-1})$.

Finally,~$\smash{V^{(13)}}$
has RM
by\/~$\bbQ(\sqrt{13})$,
$\smash{V^{(-1, \mu_7)}}$
CM
by\/~$\bbQ(\zeta_{28} + \zeta_{28}^{13}) = \bbQ(i, \zeta_7 + \zeta_7^{-1})$,
$\smash{V^{(-1, \mu_9)}}$
CM
by\/~$\bbQ(\zeta_{36} + \zeta_{36}^{17}) = \bbQ(i, \zeta_9 + \zeta_9^{-1})$,
and\/
$\smash{V^{(-1, \mu_{19})}}$
CM
by\/
$L(i)$,
for
$L \subset \bbQ(\mu_{19})$
the unique cubic~subfield.
\end{abc}
\end{cons}

\begin{rems}

\begin{iii}
\item
The~equations given actually describe singular models of the
$K3$~surfaces
to be~considered.
\item
(Evidence.)
Thanks~to Harvey's
\mbox{$p$-adic}
point counting method, we can now give a lot more numerical evidence for our conjectures than before~\cite{EJ14}.
\item
(Proven cases.)
The first three examples were already published in~\cite{EJ14}. For the
family~$\smash{V_a^{(2)}}$,
both a) and~b) were proven in \cite[Theorem~6.6]{EJ14}. The~same method applies to
$\smash{V_{{\bf a},{\bf b}}^{(-1)}}$
and provides a proof also for this~family. As~far as only~a) is concerned, using van Luijk's method one shows in each case that the geometric Picard rank is
$16$
or~$17$,
and~$16$
as soon as~b) is~true.
\item
In~particular, for each family or surface, we explicitly know 16 divisors that are linearly independent in the Picard~group. In~neither case, all of them are defined
over~$\bbQ$.
For~instance,
for~$\smash{V^{(-1, \mu_7)}}$,
$\smash{V^{(-1, \mu_9)}}$,
and
$\smash{V^{(-1, \mu_{19})}}$,
the fields of definition are
$\smash{\bbQ(\zeta_7 + \zeta_7^{-1})}$,
$\smash{\bbQ(\zeta_9 + \zeta_9^{-1})}$,
and the unique cubic subfield of
$\smash{\bbQ(\mu_{19})}$,
respectively.

On~the other hand, for
$\smash{V_a^{(2)}}$,
the field of definition of the known divisors
is~$\bbQ(\sqrt{2})$,
while it
is~$\bbQ(\sqrt{5}, \sqrt{a^2-20a-20})$
for
$\smash{V_a^{(5)}}$,
$\bbQ(\sqrt{13}, \sqrt{-3})$
for
$V^{(13)}$,
$\smash{\bbQ(\sqrt{2}, \sqrt{a}, \sqrt{b}, \sqrt{a^2-6ab+b^2})}$
for~$\smash{V_{a,b}^{(2)}}$
and
$\smash{\bbQ(\sqrt{3}, \sqrt{a^4 - 12a^3b + 30a^2b^2 - 36ab^3 + 9b^4}, \sqrt{2a^2 - 6ab + 6b^2})}$
for
$\smash{V_{a,b}^{(3)}}$.
We~do not know whether the following is more than a~coincidence.
\end{iii}
\end{rems}

\begin{eo}
All the examples listed in Conjectures~\ref{list_ex} that are supposed to have real multiplication have the property that the endomorphism field is contained in the field of definition of the Picard~group.
\end{eo}

{\em Numerical evidence. Data for the real multiplication examples.}\medskip

For each prime
$p \in \{19, \ldots, 499\}$,
we inspected all the possible specialisations of the families after reduction
modulo~$p$.
Whenever~this resulted in a non-singular surface, we computed the Weil polynomial and derived the geometric Picard~rank. This~confirmed properties \ref{arith_cons}.i) and~ii) in every~case. In~particular, only the Picard~ranks 18 and 22 are~occurring. There~is the following~statistics.

\begin{table}[H]
\small
\centerline{
\begin{tabular}{p{1.3cm}|c|c|p{1.4cm}|c|c|c}
Family & \multicolumn{6}{c}{\# relative frequency (in\,\%) of rank 22 per prime} \\[-.7mm]
       & \multicolumn{3}{c}{inert primes} & \multicolumn{3}{c}{split primes} \\[-.7mm]
       & $\min$ & average & $\max$ & $\min$ & average & $\max$ \\[-.7mm]\hline
$V_a^{(2)}$     & 0.00 &  7.42 & 25.00 & 0.00 & 6.14 & 25.00 \\
$V_a^{(5)}$     & 2.33 &  9.32 & 24.24 & 0.00 & 5.84 & 16.00 \\
$V_{a,b}^{(2)}$ & 4.75 & 12.66 & 38.89 & 1.50 & 7.39 & 22.73 \\[.85mm]
$\smash{V_{a,b}^{(3)}}$ & 0.00 &  0.00 & 0.000 & 2.33 & 7.40 & 20.59 \\
\end{tabular}}
\caption{Frequency of reduction to geometric Picard rank 22}
\end{table}

One~might want to look more closely at the primes that result in reduction to geometric Picard
rank~$22$.
In~the case of the
family~$\smash{V_a^{(2)}}$,
77 of 88 primes are~occurring. The~missing ones are
$23$,
$29$,
$31$,
$47$,
$97$,
$127$,
$193$,
$241$,
$401$,
$433$,
and~$449$.
The reductions of these primes
modulo~$8$
are
$1^6, 5, 7^4$.
The~only one being inert in the endomorphism field
$\bbQ(\sqrt{2})$
is~$29$.
For~the family
$\smash{V_{a,b}^{(2)}}$,
all 88 inspected primes do occur, while for
$\smash{V_a^{(5)}}$
87 of the 88 primes occur and only
$29$
is~missing.
However,~in the example
$\smash{V_{a,b}^{(3)}}$,
only 44 of the 88 primes result in reductions to rank 22 and in fact {\em exactly those splitting\/} in the endomorphism field
$\bbQ(\sqrt{3})$,
a phenomenon we have no explanation~for.\medskip

{\em Primes of reduction to rank 18.}
In~the case of reduction to Picard
rank~$18$,
the transcendental factor split off the Weil polynomial is of degree~four. In~our families, it was never a perfect power, cf.~\cite[Theorem~1.1]{Za93}, but always turned out to be~irreducible. It~was either the norm of a quadratic polynomial
$g \in E[t]$
over the endomorphism field or turned into a square under the operation of squaring all its~roots. The~first alternative occurred precisely at the primes split
in~$E$,
while the second option came up at the inert~primes.\bigskip

{\em Numerical evidence. Data for the complex multiplication examples.}\medskip

Let~us report only on the numerical evidence for the three isolated~examples.
These~are contained in the family
$\smash{V_{{\bf a},{\bf b}}^{(-1)}}$.
Thus,~at least CM
by~$\bbQ(\sqrt{-1})$
is~proven.
As the surfaces have geometric Picard
rank~$16$,
the transcendental part of the cohomology is of relative dimension~one. Thus,~one should expect phenomena that are very close to those known to occur for elliptic curves with complex~multiplication.

To~verify this, we run the point counting for all the primes
below~$1000$.
It~turns out that all the reductions have geometric Picard rank either
$16$
or~$22$,
where
rank~$16$
appears exactly for the primes
$1 \bmod 4$.
A~prime is ordinary if and only if it completely splits in the endomorphism~field.
Furthermore, the Frobenius eigenvalues on the transcendental part of
$\smash{H^2_\et(V_{\overline\bbF_{\!p}}, \bbQ_l)}$
are closely related to the endomorphism
field~$E$.
More~precisely,

\begin{abc}
\item
In the case that
$p$
completely splits
in~$E$,
the eigenvalues
of~$\Frob$
are contained
in~$E$.
They~are of the~form
$\pm p \frac{\overline{\pi}}{\pi}$
for~$(\pi)$
a prime
above~$p$.
Note that the endomorphism fields have class
number~$1$,
such that the ideals
above~$p$
are~principal.
\item
In the case that
$p$
splits into three primes, the characteristic polynomial
of~$\Frob$
on the transcendental lattice is
$t^6 - 3 p^2 t^4 + 3 p^4 t^2 - p^6 = (t^2 - p^2)^3$.
\item
In the case that
$p$
splits into two primes, the characteristic polynomial
of~$\Frob$
on the transcendental lattice is of the form
$t^6 + a p^2 t^3 + p^6$,
for~$a \in \bbZ$
having the property that the discriminant of
$t^2 + a t + p^2$
is minus a~square. Thus, the Weil polynomial
of~$\smash{X_{\bbF_{\!p^3}}}$
is a cube of a quadratic polynomial having its roots in
$\bbQ(i)$.
These~are of the form
$p^2 \varrho^2$,
for
$\varrho \in \bbQ(i)$
a prime
above~$p$.
Let~us note here once again that
$\bbQ(i) \subset E$.
\item
In the case that
$p$
is totally inert, the degree six transcendental factor of the Weil polynomial is equal
to~$t^6 - p^6$.
\end{abc}
The total computation took about 110 days of CPU-time on one
Intel Xeon E5-4650 processor running at 2.7\,GHz.

\end{document}